\documentclass[12pt,reqno]{amsart}
\usepackage{amsmath, amsthm, amssymb}
\usepackage[margin=1.9cm]{geometry}
\usepackage{epsfig}
\usepackage{longtable}  
\usepackage{rawfonts}
\usepackage{enumerate}
\usepackage{graphics}
\usepackage{multirow}
\usepackage{xspace}
\usepackage{graphicx}
\usepackage{pgf,tikz,pgfplots}
\usepackage{mathrsfs}
\usetikzlibrary{arrows}
\usepackage{amsmath}
\usepackage{amsfonts}
\usepackage{amssymb}
\usepackage{amsthm}
\usepackage{graphicx}
\usepackage{booktabs}
\usepackage{caption}
\usepackage{listings}
\usepackage{setspace}
\usepackage[mathscr]{eucal}
\usepackage{pgfplots}
\usepackage{hyperref}
\usepackage{wrapfig}
\usepackage{floatflt,epsfig}
\usepackage{ dsfont }
\usepackage{amscd}
\usepackage{tikz-cd}
\usepackage{fancyhdr}
\usepackage[all]{xy}
\usepackage{latexsym}
\usepackage{amscd}
\usepackage{pifont}
\usepackage{eufrak}
\usepackage{subfig}
\usepackage{easyReview}
\usepackage{listings}
\usepackage{verbatim}
\def\ZZ{{\mathbb Z}}

\newcommand{\cQ}{\mathcal{Q}}

\newcommand{\cR}{\mathcal{R}}
\newcommand{\cS}{\mathcal{S}}

\theoremstyle{plain}

\theoremstyle{theorem}
\newtheorem{defn}{Definition}[section]

\newtheorem{thm}[defn]{Theorem}

\newtheorem{coro}[defn]{Corollary}
\newtheorem{exa}[defn]{Example}

\theoremstyle{remark}

\begin{document}
	

    \title[Collections of cells, binomial ideals and combinatorics]{Collections of cells, binomial ideals and combinatorics}
    
			\author{CARMELO CISTO}
	        \address{Università di Messina, Dipartimento di Scienze Matematiche e Informatiche, Scienze Fisiche e Scienze della Terra\\
			Viale Ferdinando Stagno D'Alcontres 31\\
			98166 Messina, Italy}
			\email{carmelo.cisto@unime.it}
			
            \author{RIZWAN JAHANGIR}
			\address{Sabanci University, Faculty of Engineering and Natural Sciences, Orta Mahalle, Tuzla 34956, Istanbul, Turkey, and NUST Business School, NUST H-12 Campus, Off Srinagar Highway, Islamabad 44000, Pakistan}
			\email{rizwan@sabanciuniv.edu,
            rizwan.jahangir@nbs.nust.edu.pk}
   
			\author{FRANCESCO NAVARRA}
			\address{Sabanci University, Faculty of Engineering and Natural Sciences, Orta Mahalle, Tuzla 34956, Istanbul, Turkey}
			\email{francesco.navarra@sabanciuniv.edu}

	

	\begin{abstract}  
    In this paper we provide a description of the package \textit{PolyominoIdeals} for \textit{Macaulay2} that allows to deal with collections of cells, polyominoes and related binomial ideals.
	\end{abstract}

	\maketitle
	
\section*{Introduction}
Collections of unit squares in the integer lattice give rise to a broad and versatile family of combinatorial objects. By considering collections of such cells in the plane, one can construct a wide variety of geometric figures with different shapes whose roles exhibit deep connections in commutative algebra, combinatorics, and computational algebra.

A particularly well-studied class of collections of cells is that of \emph{polyominoes}, which are finite collections of unit squares joined edge to edge. 
The term “polyomino” was coined by Solomon W. Golomb in 1953, and these configurations have since been investigated mainly in combinatorial mathematics, especially in relation to tiling problems of the plane (see \cite{golomb}). 
Beyond their combinatorial significance, polyominoes have also found applications in several other fields, including theoretical computer science, statistical physics, and discrete geometry (see \cite{GG, GG2}).

More recently, since 2012, polyominoes and, more generally, collections of cells have been studied from an algebraic–combinatorial perspective. In \cite{Q}, Ayesha Asloob Qureshi established a connection between combinatorial commutative algebra and collections of cells by associating to each collection the binomial ideal generated by its inner $2$-minors. If the generators are restricted to those corresponding only to adjacent cells, one obtains the so-called \emph{adjacent $2$-minor ideals}, introduced by Herzog and Hibi in \cite{HH2012}. 

For a more comprehensive study of binomial ideals, the reader may consult \cite{ES, binomial ideals} and the \textit{Macaulay2} package \textit{Binomials}, developed by T.~Kahle (\cite{Binomials}).

   In this paper, we present a detailed description of the package \textit{PolyominoIdeals} (\cite{Polyo}) for the computer algebra system \textit{Macaulay2} (\cite{M2}). More precisely, once a collection of cells is encoded as a list of the lower-left corners of its cells, the package provides a variety of combinatorial and algebraic functions, which are described in Sections~1 and~2, respectively. Finally, in Section~3, we offer an overview of rook theory together with several related functions.

	\section{Collections of cells, polyominoes and some combinatorial functions}\label{Section: Introduction}

This section formally introduces the notion of collections of cells, polyominoes, and their main combinatorial properties, together with the associated structural–geometric functions implemented in \textit{PolyominoIdeals} package.

   \subsection*{Collections of cells.} Let $(i,j),(k,l)\in \mathbb{Z}^2$. We say that $(i,j)\leq (k,l)$ if $i\leq k$ and $j\leq l$. 
Consider $a=(i,j)$ and $b=(k,l)$ in $\mathbb{Z}^2$ with $a\leq b$. The set 
\[
[a,b]=\{(m,n)\in \mathbb{Z}^2 : i\leq m\leq k,\ j\leq n\leq l\}
\]
is called an \textit{interval} of $\mathbb{Z}^2$. 
If $i<k$ and $j<l$, then $[a,b]$ is a \textit{proper} interval. In this case, we call $a$ and $b$ the \textit{diagonal corners} of $[a,b]$, and $c=(i,l)$ and $d=(k,j)$ the \textit{anti-diagonal corners}. 
If $j=l$ (or $i=k$), then $a$ and $b$ are said to be in \textit{horizontal} (or \textit{vertical}) \textit{position}. 
We denote by $]a,b[$ the set $\{(m,n)\in \mathbb{Z}^2 : i<m<k,\ j<n<l\}$.

A proper interval $C=[a,b]$ with $b=a+(1,1)$ is called a \textit{cell} of $\mathbb{Z}^2$. 
The points $a$, $b$, $c$, and $d$ are called the \textit{lower left}, \textit{upper right}, \textit{upper left}, and \textit{lower right corners} of $C$, respectively. 
The sets $\{a,c\}$, $\{c,b\}$, $\{b,d\}$, and $\{a,d\}$ are the \textit{edges} of $C$. 
We set $V(C)=\{a,b,c,d\}$ and $E(C)=\{\{a,c\},\{c,b\},\{b,d\},\{a,d\}\}$.

 A cell in the plane $\mathbb{Z}^2$ is uniquely identified by its lower-left corner; hence a collection of cells can be encoded as a list of these lower-left corners.

\begin{figure}[h!]
    \centering
    \includegraphics[scale=0.9]{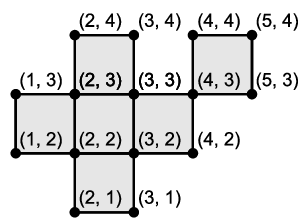}
    \caption{A collection of cells $\cQ$.}
    \label{Figure: Polyomino}
\end{figure}

For instance, the collection of cells in Figure~\ref{Figure: Polyomino} is represented by
\begin{center}
\texttt{L = \{\{1, 2\}, \{2, 1\}, \{2, 2\}, \{2, 3\}, \{3, 2\}, \{4,3\}\}}.
\end{center}

The command \texttt{cellCollection L} creates an object of type \texttt{CollectionOfCells}.

\subsection*{Vertices, edges and rank.} If $\mathcal{Q}$ is a non-empty collection of cells in $\mathbb{Z}^2$, we define the sets of vertices and edges of $\mathcal{Q}$ by 
\[
V(\mathcal{Q})=\bigcup_{C\in \mathcal{Q}}V(C), \qquad 
E(\mathcal{Q})=\bigcup_{C\in \mathcal{Q}}E(C),
\]
and $\mathrm{rank}(\mathcal{Q})$ denotes the number of cells in $\mathcal{Q}$. The functions \texttt{polyoVertices Q} and \texttt{polyoEdges Q} compute the list of the vertices and edges of \texttt{Q} and \texttt{rankCollection Q} determines the rank of \texttt{Q}.

\subsection*{Maximal edge intervals.}
If $\cQ$ is a collection of cells, an interval $[a,b]$ with $a=(i,j)$, $b=(k,j)$, and $i<k$ is called a \textit{horizontal edge interval} of $\cQ$ if the set $\{(\ell,j),(\ell+1,j)\}$ is an edge of a cell of $\cQ$, for all $\ell=i,\dots,k-1$. 
If $\{(i-1,j),(i,j)\}$ and $\{(k,j),(k+1,j)\}$ do not belong to $E(\cQ)$, then $[a,b]$ is a \textit{maximal horizontal edge interval} of $\cQ$. The vertical and maximal vertical edge intervals are defined analogously.

The functions \texttt{maximalVerticalEdgeIntervals Q} and \texttt{maximalHorizontalEdgeIntervals Q} compute the maximal vertical and horizontal edge interval of \texttt{Q}, respectively.

\subsection*{Connectedness.} 
\begin{figure}[h!]
\centering
\subfloat{\includegraphics[scale=0.6]{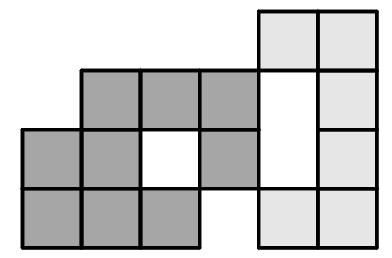}}\qquad\qquad\qquad
\subfloat{\includegraphics[scale=0.6]{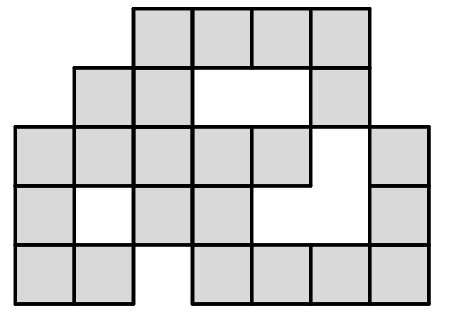}}
\caption{A weakly connected collection of cells with two connected components on the left and a polyomino on the right.}
\label{Figura: A polyomino+weakly connected}
\end{figure}

Let $\cQ$ be a collection of cells. If $C$ and $D$ are two distinct cells of $\mathcal{Q}$, a \textit{walk} from $C$ to $D$ in $\cQ$ is a sequence $\mathcal{C}: C=C_1,\dots,C_m=D$ of cells of $\mathbb{Z}^2$ such that $C_i\cap C_{i+1}$ is an edge of both $C_i$ and $C_{i+1}$ for all $i=1,\dots,m-1$. 
If $C_i\neq C_j$ for all $i\neq j$, then $\mathcal{C}$ is called a \textit{path} from $C$ to $D$. 
Two cells $C$ and $D$ are \textit{connected} in $\cQ$ if there exists a path of cells in $\cQ$ joining them. A collection of cells $\cQ$ is called a \textit{polyomino} if every couple of cells of $\cQ$ is connected in $\cQ$. In this case, we also say that $\cQ$ is a \textit{connected} collection of cells.\\
The collection $\cQ$ is said to be \textit{weakly connected} if, for any two cells $C$ and $D$ in $\cQ$, there exists a sequence $\mathcal{C}: C=C_1,\dots,C_m=D$ of cells in $\cQ$ such that $V(C_i)\cap V(C_{i+1})\neq \emptyset$ for all $i=1,\dots,m-1$. 
A subset $\cQ'\subseteq \cQ$ is called a \textit{connected component} of $\cQ$ if $\cQ'$ is a polyomino and it is maximal with respect to inclusion; that is, if $A\in \cQ\setminus \cQ'$, then $\cQ'\cup \{A\}$ is not a polyomino. 
Clearly, every polyomino is a weakly connected collection of cells. 
For instance, see Figure~\ref{Figura: A polyomino+weakly connected}.

The function \texttt{collectionIsConnected Q} tests whether \texttt{Q} is a polyomino and, if not, provides the number of connected components. Moreover, the function \texttt{connectedComponentsCells Q} computes the connected components of \texttt{Q}.

To determine whether a collection of cells $\cQ$ is a polyomino, it suffices to consider the graph $G_\cQ$ associated with $\cQ=\{C_1,\dots,C_n\}$, where 
$V(G_\cQ)=\{1,\dots,n\}$ and 
$E(G_\cQ)=\{\{i,j\} : i\neq j,\ C_i \text{ share an edge with } C_j\}$. 
The connectedness of $\mathcal{Q}$ is equivalent to the connectedness of the graph $G_{\mathcal{Q}}$.   The function \texttt{cellGraph Q} constructs the graph associated with \texttt{Q}, as defined above, by making use of the \textit{Graphs} package (see \cite{Graphs}).

\subsection*{Inner intervals.} Let $A$ and $B$ be two cells in $\mathbb{Z}^2$ with lower-left corners $a=(i,j)$ and $b=(k,l)$, respectively, and assume $a\leq b$. 
The \textit{cell interval} $[A,B]$ is the set of all cells in $\mathbb{Z}^2$ whose lower-left corner $(r,s)$ satisfies $i\leq r\leq k$ and $j\leq s\leq l$. 
If $(i,j)$ and $(k,l)$ lie in horizontal (respectively, vertical) position, then $A$ and $B$ are said to be in \textit{horizontal} (respectively, \textit{vertical}) position.

Observe that any interval of $\mathbb{Z}^2$ uniquely determines a cell interval in $\mathbb{Z}^2$, and conversely; thus, to every interval $I$ of $\mathbb{Z}^2$ we associate the corresponding cell interval $\cQ_I$.

If $\cQ$ is a collection of cells, a proper interval $[a,b]$ of $\ZZ^2$ is an \emph{inner interval} of $\cQ$ if every cell of $\cQ_{[a,b]}$ is contained in $\cQ$. The function \texttt{innerInterval(\{a,b\},\{c,d\},Q)} determines if $[(a,b),(c,d)]$ is an inner interval of \texttt{Q}, where \{\texttt{a}$, $\texttt{b}\}, \{\texttt{c}, \texttt{d}\} encode respectively the lattice points $(a,b)$ and $(c,d)$.

\subsection*{Convexity.} A collection of cells $\cQ$ is \textit{row convex} if, for any two distinct cells $A$ and $B$ of $\mathcal{Q}$ in a horizontal position, the interval $[A,B]$ is contained in $\mathcal{Q}$. 
Equivalently, for any two distinct vertices $a$ and $b$ of $\mathcal{Q}$ in horizontal position and with $a<b$, the interval $[a, b+(1,1)]$ is an inner interval of $\mathcal{Q}$. 
This is the criterion used to implement the convexity test. 
Similarly, $\mathcal{Q}$ is \textit{column convex} if the same property holds for cells in vertical position. 
The collection $\mathcal{Q}$ is \textit{convex} if it is both row and column convex.

The functions \texttt{isRowConvex Q} and \texttt{isColumnConvex Q} check the row and column convexity of \texttt{Q}, respectively. The function \texttt{isConvex Q} tests whether \texttt{Q} is convex.

\subsection*{Simplicity.} A collection of cells $\mathcal{Q}$ is called \textit{simple} if, for any two cells $C$ and $D$ not in $\cQ$, there exists a path of cells not belonging to $\cQ$ connecting $C$ and $D$. 
A finite collection of cells $\mathcal{H}$ not contained in $\cQ$ is a \textit{hole} of $\cQ$ if any two cells of $\mathcal{H}$ are connected within $\mathcal{H}$ and $\mathcal{H}$ is maximal with respect to inclusion. 
For example, the polyomino in Figure~\ref{Figura: A polyomino+weakly connected}(A) is not simple and contains three holes. Each hole of $\cQ$ is itself a simple polyomino, and $\cQ$ is simple if and only if it has no holes. To test the simplicity of $\cQ$, one may proceed into two steps. 
Firstly, let $[a,b]$ be an interval such that $V(\cQ)\subset [a,b]$, and define $\cR=\cQ_I$ where $I=[a-(1,1), b+(1,1)]$. Secondly, then we check whether the set $\cR \setminus \cQ$ is connected.

The function \texttt{collectionIsSimple Q} determines whether \texttt{Q} is simple.

\subsection*{Random collections of cells.} Sometimes it is useful to generate random collections of cells or random polyominoes. 
In this setting, functions that determine properties such as simplicity or convexity become particularly valuable: 
for large collections, these properties can be assessed without the need to visualize the configuration. 
To construct a random polyomino of rank $n$, we employ the function \texttt{random m}, which gives a random integer between $0$ and $m$.
More precisely, we start from the initial collection \texttt{Q = \{\{1,1\}\}} and, while the rank of \texttt{Q} is smaller than $n$, 
we consider the set of cells in $\ZZ^2$ that are not already in \texttt{Q} but share an edge with a cell of \texttt{Q}. 
This set, denoted \texttt{available} in the code, is never empty. 
We then select a random cell from \texttt{available} and include it in \texttt{Q}. 
A similar procedure can be used to generate a random collection of cells that is not necessarily a polyomino, where in \texttt{available} is the list of cells in $\ZZ^2$ that are not already in \texttt{Q} but share a vertex with a cell of \texttt{Q}. 

 The functions \texttt{randomCollectionWithFixedRank n} and \texttt{randomPolyominoWithFixedRank n} generates a random collection of cells or polyomino with rank $n$, respectively.
 
For random collections of cells or random polyominoes with rank bounded by $n$, one may first choose a random integer between $1$ and $n$ and then apply the same construction using this value as the target rank. The functions \texttt{randomCollectionOfCells n} and \texttt{randomPolyomino n} generates a random collection of cells and polyomino with rank at most $n$, respectively.

\subsection*{An example.} We now conclude the section by showing an explicit example illustrating the usage of the described functions. 

\begin{footnotesize}
\begin{lstlisting}
i1 : Q = randomCollectionOfCells 10 
o1 = {{1, 1}, {0, 0}, {0, 2}, {2, 2}, {1, 0}, {0, 3}}
o1 : CollectionOfCells
i2 : polyoVertices Q
o2 = {(3, 2), (3, 3), (2, 0), (2, 1), (2, 2), (2, 3), (1, 0), (1, 1),
      (1, 2), (1, 3), (1, 4), (0, 0), (0, 1), (0, 2), (0, 3), (0, 4)}
i3 : innerInterval({1,1},{3,4},Q)
o3 = false
i5 : collectionIsConnected Q
o5 = (false, 3)        -- Q is not a polyomino and it has 3 connected components 
i6 : connectedComponentsCells Q
o6 = {{{1, 1}, {1, 0}, {0, 0}}, {{0, 2}, {0, 3}}, {{2, 2}}
i7 : isRowConvex Q
o7 = false
i8 : isColumnConvex Q
o8 = false
i9 : collectionIsSimple Q
o9 = true
i10 : rankCollection Q
o10 = 6
\end{lstlisting}
\end{footnotesize}


\section{Commutative algebra related to collections of cells}\label{Section: Functions}

In this section, we describe the algebraic functions included in the package, togheter the related options.

\subsection*{Inner 2-minor ideal.} We recall the definition of the \emph{inner $2$-minor ideal} of a collection of cells, following~\cite{Q}.  

Let $\cQ$ be a collection of cells. Let $S_{\cQ}=K[x_v \mid v \in V(\cQ)]$, where $K$ is a field. If $[a,b]$ is an inner interval of $\cQ$, and $a,b$ and $c,d$ are its diagonal and anti-diagonal corners respectively, then the binomial $x_a x_b - x_c x_d$ is called an \emph{inner $2$-minor} of $\cQ$. The ideal $I_{\cQ}$ generated by all inner $2$-minors of $\cQ$ is the \emph{inner $2$-minor ideal} of $\cQ$. When $\cQ$ is a polyomino, $I_{\cQ}$ is also called the \emph{polyomino ideal} of $\cQ$.

For the polyomino $\cQ$ in Figure~\ref{Figure: Polyomino}, the ideal $I_{\cQ}$ is
\[
\begin{aligned}
I_\cQ = (&x_{3,3}x_{2,2} - x_{3,2}x_{2,3},\ 
x_{3,4}x_{2,1} - x_{3,1}x_{2,4},\ 
x_{4,3}x_{1,2} - x_{4,2}x_{1,3},\ 
x_{4,3}x_{2,2} - x_{4,2}x_{2,3},\\
&x_{3,2}x_{2,1} - x_{3,1}x_{2,2},\ 
x_{3,4}x_{2,2} - x_{3,2}x_{2,4},\ 
x_{2,3}x_{1,2} - x_{2,2}x_{1,3},\ 
x_{4,3}x_{3,2} - x_{4,2}x_{3,3},\\ 
&x_{3,3}x_{2,1} - x_{3,1}x_{2,3},\ 
x_{3,4}x_{2,3} - x_{3,3}x_{2,4},\ 
x_{3,3}x_{1,2} - x_{3,2}x_{1,3},\
x_{4,3}x_{5,4} - x_{4,4}x_{5,3}).
\end{aligned}
\]

The function \texttt{polyoIdeal} returns the inner $2$-minor ideal $I_{\cQ}$.

\subsection*{Adjacent 2-minor ideal.}

We now recall the \emph{adjacent $2$-minor ideal}, following~\cite{HH2012}.  
Let $\cQ$ be a collection of cells, and let $S_{\cQ}=K[x_v \mid v \in V(\cQ)]$ as above.

The ideal $I_{\mathrm{adj}}(\cQ)$ is generated by all binomials $x_a x_b - x_c x_d$ such that $[a,b]$ consists of a single cell of $\cQ$, and $c,d$ are the anti-diagonal corners of $[a,b]$.  
This ideal is called the \emph{adjacent $2$-minor ideal} of $\cQ$. The quotient $K[\cQ]_{\mathrm{adj}} = S_{\cQ} / I_{\mathrm{adj}}(\cQ)$ is the corresponding \emph{adjacent coordinate ring} of $\cQ$.

For the polyomino $\cQ$ in Figure~\ref{Figure: Polyomino}, the ideal $I_{\mathrm{adj}}(\cQ)$ is
\[
\begin{aligned}
I_{\mathrm{adj}}(\cQ) = (&x_{3,2}x_{2,1} - x_{3,1}x_{2,2},\ 
x_{2,3}x_{1,2} - x_{2,2}x_{1,3},\ 
x_{3,3}x_{2,2} - x_{3,2}x_{2,3},\\
&x_{4,3}x_{3,2} - x_{4,2}x_{3,3},\ 
x_{3,4}x_{2,3} - x_{3,3}x_{2,4},\
x_{4,3}x_{5,4} - x_{4,4}x_{5,3}).
\end{aligned}
\]

The function \texttt{adjacent2MinorIdeal} returns the adjacent $2$-minor ideal $I_{\mathrm{adj}}(\cQ)$.

\subsection*{Collection of cells as a matrix.}
Let $\cQ$ be a collection of cells and $[(a,b),(c,d)]$ the smallest interval of $\ZZ^2$ containing $\cQ$. The matrix $M(\cQ)$ has $d-b+1$ rows and $c-a+1$ columns, with $M(\cQ)_{i,j}=x_{(i,j)}$ if $(i,j)$ is a vertex of $\cQ$, and $0$ otherwise.

Consider the same polyomino shown in Figure~\ref{Figure: Polyomino}. The associated matrix is produced using the {\tt polyoMatrix} function as follows:
\begin{footnotesize}
	\begin{lstlisting}
i1 : Q = {{1, 2}, {2, 1}, {2, 2}, {2, 3}, {3, 2}};
i2 : polyoMatrix(Q);
o3 : | 0       x_(2,4) x_(3,4) 0       |
     | x_(1,3) x_(2,3) x_(3,3) x_(4,3) |
     | x_(1,2) x_(2,2) x_(3,2) x_(4,2) |
     | 0       x_(2,1) x_(3,1) 0       |
    \end{lstlisting}
    \end{footnotesize}

\noindent This function is primarly used for implementing the function \texttt{polyoRingReduced}, which is essential for the option used when \texttt{RingChoice} takes a value other than $1$ (see Subsection~\ref{subsection:RingChoice}).

\subsection*{Toric representations.}

\noindent Let $\cQ$ be a weakly connected collection of cells. We introduce a toric ideal naturally associated with $\cQ$, extending the construction given in \cite{MRR} for polyominoes. Consider the following total order on $V(\cQ)$: for $a=(i,j)$ and $b=(k,l)$, set $a\succ b$ if either $i>k$, or $i=k$ and $j>l$. If $\mathcal{H}$ is a hole of $\cQ$ and $e\in V(\cQ)$, we say that $e$ is the \textit{lower-left corner} of $\mathcal{H}$ is $e$ is the minimum, with respect to $\prec$, of the vertices of $\mathcal{H}$. Let $\mathcal{H}_1,\dots,\mathcal{H}_r$ be the holes of $\cQ$, and let $e_k=(i_k,j_k)$ denote the lower-left corner of $\mathcal{H}_k$. For $k\in K=[r]$, define
\[
F_k=\{(i,j)\in V(\cQ) : i\le i_k,\ j\le j_k\}.
\]

Let $\{V_i\}_{i\in I}$ be the set of all maximal vertical edge intervals of $\cQ$, and let $\{H_j\}_{j\in J}$ be the set of all maximal horizontal edge intervals of $\cQ$. Let $\{v_i\}_{i\in I}$, $\{h_j\}_{j\in J}$, and $\{w_k\}_{k\in K}$ be three sets of variables. Consider the map
\[
\alpha : V(\cQ) \longrightarrow K[h_i, v_j, w_k : i\in I,\ j\in J,\ k\in K]
\]
\[
a \longmapsto h_i v_j \prod_{k\in K} w_k^{\epsilon_k}.
\]
where $a\in H_i\cap V_j$ and, for every $k\in K$, we set $\epsilon_k=1$ if $a\in F_k$ and $\epsilon_k=0$ if $a\notin F_k$. The toric ring $T_{\cQ}$ associated with $\cQ$ is defined by
\[
T_{\cQ} = K[\alpha(a) : a\in V(\cQ)].
\]
The homomorphism $\psi : S_\cQ \rightarrow T_{\cQ}$ given by $x_a \mapsto \alpha(a)$ is surjective, and the toric ideal $J_{\cQ}$ is its kernel. This construction naturally extends the one presented in \cite{CNU1} and \cite{QSS}.

\begin{thm}\cite[Theorem 3.3]{CNU1}\label{Theorem: primality and polyotoric}
Let $\cQ$ be a simple and weakly connected collection of cells. Then $I_{\cQ}=J_{\cQ}$.
\end{thm}

 The function {\tt PolyoToric(Q,H)} computes the toric ideal $J_{\cQ}$ defined above, where {\tt Q} is the list encoding the collection of cells and {\tt H} is the list of the lower-left corners of the holes. We illustrate this with two examples below.

\begin{exa}\rm
Consider the simple and weakly connected collection $\cQ$ of cells in Figure~\ref{Figura: polyotoric} (A). We compute the ideal $I_{\cQ}$ using {\tt polyoIdeal(Q)}, the toric ideal $J_{\cQ}$ with {\tt polyoToric(Q,\{\})}, and then compare the two ideals. To verify equality, we must first bring {\tt J=polyoToric(Q,\{\})} into the ring {\tt R} of {\tt polyoIdeal(Q)} using the command {\tt substitute(J,R)}. In accordance with Theorem~\ref{Theorem: primality and polyotoric}, we obtain $I_{\cQ}=J_{\cQ}$.

\begin{figure}[h!]
\centering
\subfloat[]{\includegraphics[scale=0.5]{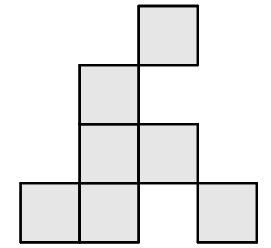}}\qquad\qquad\qquad
\subfloat[]{\includegraphics[scale=0.5]{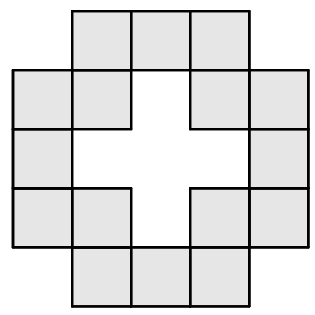}}
\caption{Some collections of cells.}
\label{Figura: polyotoric}
\end{figure}

\begin{footnotesize}
\begin{lstlisting}
i1 : Q={{1, 1}, {2, 2}, {2, 1}, {3, 2}, {2, 3}, {4, 1}, {3, 4}};
i2 : I=polyoIdeal(Q);
i3 : J=polyoToric(Q,{});
i4 : R=ring I;
i5 : J=substitute(J,R);
o5 : Ideal of R
i6 : J==I
o6 = true
\end{lstlisting}
\end{footnotesize}
\end{exa}

\begin{exa}\rm
\noindent Consider now the polyomino $\cQ$ in Figure~\ref{Figura: polyotoric} (B) . The polyomino ideal is not prime (see \cite{CN}), hence $I_{\cQ} \subset J_{\cQ}$ since $I_{\cQ}=(J_{\cQ})_2$ (Lemma~3.1, \cite{MRR}). We can also compute the set of binomials generating $J_{\cQ}$ but not $I_{\cQ}$.

\begin{footnotesize}
\begin{lstlisting}
i1 : Q={{2, 1}, {2, 2}, {1, 2}, {1, 3}, {1, 4}, {2, 4}, {2, 5}, {3, 5}, 
        {4, 5}, {4, 4}, {5, 4}, {5, 3}, {5, 2}, {4, 2}, {4, 1}, {3, 1}};
i2 : I=polyoIdeal(Q);
i3 : J=polyoToric(Q,{{2,3}});
i4 : R=ring I;
i5 : J=substitute(J,R);
i6 : J==I
o6 = false
i7 : select(first entries mingens J,f->first degree f>=3)
o7 = {x   x   x   x    - x   x   x   x   }
       6,5 5,1 2,6 1,2    6,2 5,6 2,1 1,5
\end{lstlisting}
\end{footnotesize}
\end{exa}

\subsection*{Lattice ideal.} 

Let $\cQ$ be a collection of cells. For every $a \in V(\cQ)$ we denote by $\mathbf{v}_a$ the vector in $\mathbb{Z}^{|V(\cQ)|}$ whose $a$-th coordinate is $1$ and all other coordinates are $0$.  
If $[a,b]\in \cQ$ is an interval with diagonal corners $a,b$ and anti-diagonal corners $c,d$, we set  
\[
\mathbf{v}_{[a,b]}=\mathbf{v}_a+\mathbf{v}_b-\mathbf{v}_c-\mathbf{v}_d \in \mathbb{Z}^{|V(\cQ)|}.
\]
We define $\Lambda_\cQ$ as the sublattice of $\mathbb{Z}^{|V(\cQ)|}$ generated by the vectors $\mathbf{v}_I$ for all cells $I\in \cQ$. 

Let $n=|V(\cQ)|$. For $\mathbf{v}\in \mathbb{N}^n$ we use the standard notation $x^\mathbf{v}$ to denote the monomial of $S_\cQ$ with exponent vector $\mathbf{v}$.  
Given a vector $\mathbf{e}\in \mathbb{Z}^n$, we denote by $\mathbf{e}^+$ the vector obtained by replacing all negative entries of $\mathbf{e}$ with $0$, and set $\mathbf{e}^- = -(\mathbf{e}-\mathbf{e}^+)$. 

Let $L_\cQ$ be the lattice ideal associated with $\Lambda_\cQ$, namely the binomial ideal in $S_\cQ$ given by
\[
L_\cQ = \bigl(\, x^{\mathbf{e}^+} - x^{\mathbf{e}^-} \ \mid\ \mathbf{e}\in \Lambda_\cQ \,\bigr).
\]

In studying the primality of $I_{\cQ}$, two results are important to mention.

\begin{thm}\cite[pp.288]{Q}
    Let $\cQ$ be a collection of cells $\cQ$. Then $L_\cQ$ is prime. 
\end{thm}

\begin{coro}\cite[Theorem~3.6]{Q}
    Let $\cQ$ be a collection of cells $\cQ$. Then, $I_\cQ$ is a prime ideal if and only if $I_\cQ = L_\cQ$.
\end{coro}

The function \texttt{polyoLattice} computes the lattice ideal $L_{\mathcal{Q}}$ associated with $\mathcal{Q}$, making use of the \textit{gfanInterface} package (\cite{Lattice}).

\begin{exa}\rm
Here we present an example involving two collections of cells, $\cQ_1$ and $\cQ_2$, encoded respectively by \texttt{\{\{2,1\}, \{1,2\}, \{3,2\}, \{2,3\}\}} and \texttt{\{\{2,1\}, \{1,2\}, \{3,2\}, \{2,3\}, \{2,2\}\}}. The inner $2$-minors ideal of $\cQ_1$ is not prime, whereas that of $\cQ_2$ is prime. Consequently, $I_{\cQ_1}\subsetneq L_{\cQ_1}$ while $I_{\cQ_2}=L_{\cQ_2}$.

\begin{footnotesize}
\begin{lstlisting}
i1 : Q1 = cellCollection {{2,1}, {1,2}, {3,2}, {2,3}};
i2 : J1 = polyoLattice Q1; 
i3 : first entries mingens J1
o3 = {x   x   -x   x   , x   x   -x   x   , x   x   -x   x   , 
       3,3 2,4  3,4 2,3   3,2 4,3  3,3 4,2   3,1 2,2  3,2 2,1

      x   x   -x   x   , x   x   x   x   -x   x   x   x   }
       1,2 2,3 1,3 2,2    1,2 3,1 2,4 4,3   1,3 3,4 2,1 4,2

i4 : Q2 = cellCollection {{2,1}, {1,2}, {3,2}, {2,3}, {2,2}};
i5 : J2 = polyoLattice Q2;
i6 : first entries mingens J2
o6 = {x   x   -x   x   , x   x   -x   x   , x   x   -x   x   ,
       2,2 4,3  2,3 4,2   3,3 2,4  3,4 2,3   3,2 4,3  3,3 4,2

      x   x   -x   x   , x   x   -x   x   , x   x   -x   x   ,
       3,2 2,4  3,4 2,2   3,2 2,3  3,3 2,2   3,1 2,4  3,4 2,1

      x   x   -x   x   , x   x   -x   x   , x   x   -x   x   ,
       3,1 2,3  3,3 2,1   3,1 2,2  3,2 2,1   1,2 4,3  1,3 4,2

      x   x   -x   x   , x   x   -x   x   }
       1,2 2,3  1,3 2,2   1,2 3,3  1,3 3,2
\end{lstlisting}
\end{footnotesize}
\end{exa}

\subsection*{Options: \texttt{Field}, \texttt{TermOrder} and \texttt{RingChoice}.} \label{subsection:RingChoice}


We describe three options of the two functions \texttt{polyoIdeal} and \texttt{adjacent2MinorIdeal}.

Let $\cQ$ be a collection of cells. The option \texttt{Field} for the function \texttt{polyIdeal} allows one to change the base ring of the polynomial ring embedded in $I_\cQ$. Any base ring provided by \textit{Macaulay2} (\cite{M2}) can be selected. The option \texttt{TermOrder} allows modifying the monomial order of the ambient ring of $I_\cQ$ as specified by the function \texttt{polyoIdeal}. By default, the lexicographic order is used, but one may replace it with other monomial orders available in \textit{Macaulay2} (\cite{M2}). 

The option \texttt{RingChoice} allows the user to choose between two possible ambient rings for $I_{\cQ}$.\\
If \texttt{RingChoice} is set to 1 (which is also the default), the function \texttt{polyoIdeal} returns the ideal $I_{\cQ}$ inside the polynomial ring $S_{\cQ} = K[x_a : a \in V(\cQ)]$, where $K$ is a field and the monomial order is determined by the \texttt{TermOrder} option. This order is induced by the following ordering of the variables: $x_a > x_b$ for $a = (i,j)$ and $b = (k,l)$ if $i > k$, or if $i = k$ and $j > l$.\\
We now describe the ambient ring when \texttt{RingChoice} is assigned a value different from~1. Consider the edge ring $R = K[s_i t_j : (i,j) \in V(\cQ)]$ associated with the bipartite graph $G$ having vertex set $\{s_1,\dots,s_m\} \cup \{t_1,\dots,t_n\}$, where each vertex $(i,j) \in V(\cQ)$ corresponds to the edge $\{s_i, t_j\}$ of $G$. Let $S = K[x_a : a \in V(\cQ)]$ and let $\phi : S \rightarrow R$ be the $K$-algebra homomorphism defined by $\phi(x_{ij}) = s_i t_j$ for all $(i,j) \in V(\cQ)$, and set $J_{\cQ} = \ker(\phi)$. By Theorem~2.1 of \cite{Q}, we have $I_{\cQ} = J_{\cQ}$ whenever $\cQ$ is weakly connected and convex. In this situation, \cite{HO} shows that the generators of $I_{\cQ}$ form the reduced Gr\"obner basis with respect to a suitable order $<$, and in particular the initial ideal $\mathrm{in}_< (I_{\cQ})$ is squarefree and generated in degree two. Following the argument in \cite{HO}, the implemented routine constructs the polynomial ring $S_{\cQ}$ equipped with the monomial order $<$. The following example illustrates this.

    \begin{footnotesize}
\begin{lstlisting}
i1 : Q = cellCollection {{1,3},{2,3},{2,4},{3,4},{2,5},{3,3},{3,2},{4,4} }
o1 = {{1, 3}, {2, 3}, {2, 4}, {3, 4}, {2, 5}, {3, 3}, {3, 2}, {4, 4}}
i2 : I = polyoIdeal(Q, RingChoice=>2);
i3 : gb I
o3 = GroebnerBasis[status: done; S-pairs encountered up to degree 2]
i4 : I = polyoIdeal(Q, RingChoice=>1);
i5 : gb I
o5 = GroebnerBasis[status: done; S-pairs encountered up to degree 3]
\end{lstlisting}
\end{footnotesize}

	\section{Rook polynomial theory and a connection with the Hilbert-Poincar\'e series}

The study of rook theory has long been an active and appealing area of research. The \textit{rook problem} concerns determining the number of ways to place $k$ non-attacking rooks on a polyomino $\cQ$, and remains widely open. The placement of non-attacking rooks on a skew diagram corresponds to the enumeration of permutations with certain restrictions; this idea was introduced by Kaplansky and Riordan~\cite{KR} and later developed by Riordan~\cite{R}. For a comprehensive treatment of permutations with forbidden positions, we refer the reader to Stanley~\cite[Chapter~2]{S1}.

\subsection*{Standard rook polynomial.} We introduce the standard rook polynomial, and the related functions. Let us begin by introducing some definitions.

\begin{figure}[h]
    \centering
    \subfloat[]{\includegraphics[scale=0.8]{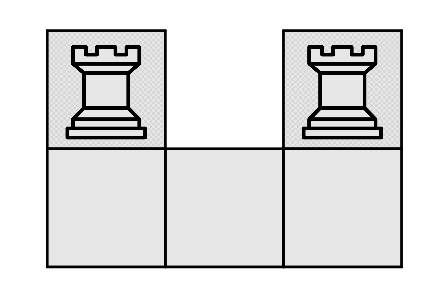}}\quad\qquad
    \subfloat[]{\includegraphics[scale=0.8]{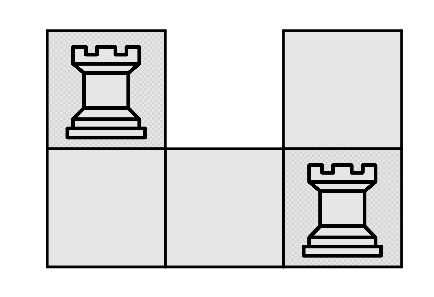}}
    \caption{Positions of two rooks in a polyomino.}
    \label{Figure: exa attacking rooks}
\end{figure}

Two rooks $R_1$ and $R_2$ are said to be in \textit{standard attacking position}, or are \textit{standard attacking rooks} in $\cQ$, if there exist two cells of $\cQ$ in horizontal or vertical position containing $R_1$ and $R_2$. Conversely, two rooks are in \textit{non-attacking position}, or are \textit{non-attacking rooks} in $\cQ$, if they are not in attacking position. For instance, examples of two standard attacking rooks and two standard non-attacking rooks are illustrated in Figures~\ref{Figure: exa attacking rooks}(A) and (B), respectively.

The \textit{standard rook number} $r_{\mathrm{st}}(\cQ)$ is the maximum number of rooks that can be placed in $\cQ$ in standard non-attacking positions. We denote by $\mathcal{R}_{\mathrm{st}}(\cQ,k)$ the set of all configurations of $k$ rooks in standard non-attacking position in $\cQ$, and set $(r_{\mathrm{st}})_k = \lvert \mathcal{R}_{\mathrm{st}}(\cQ,k) \rvert$ for all $k \in \{0, \dots, r_{\mathrm{st}}(\cQ)\}$ (with the convention $(r_{\mathrm{st}})_0 = 1$). The \textit{standard rook polynomial} of $\cQ$ is the polynomial in $\mathbb{Z}_{> 0}[t]$ defined by
$$
r_{\mathrm{st}}(\cQ)(t) = \sum_{k=0}^{r_{\mathrm{st}}(\cQ)} (r_{\mathrm{st}})_k t^k.
$$
For instance, the polyomino in Figure~\ref{Figure: exa attacking rooks} has $r_{\mathrm{st}}(\cQ) = 2$ and $r_{\mathrm{st}}(\cQ)(t) = 4 t^2 + 5t + 1$.

  The functions \texttt{standardRookNumber Q} and \texttt{standardRookPolynomial Q} computes the standard rook number and the standard rook polynomial of \texttt{Q}, respectively. Moreover, \texttt{standardNonAttackingRookConfigurations Q} returns all standard non-attacking rook configurations of \texttt{Q}.

\subsection*{A variant of the rook polynomial and Commutative Algebra.} A recent line of research has established a novel connection showing that the Hilbert--Poincar\'e series of $K[\cQ]$ is closely related to the rook polynomial and to one of its variants.

For the sake of completeness, we recall the definition of the Hilbert--Poincar\'e series of a graded ideal.
  
Let $R = K[x_1, \dots, x_n]$ be a polynomial ring over a field $K$, and let $I$ be a homogeneous ideal of $R$. By the classical Hilbert Syzygy Theorem, it is well known that $I$ admits a minimal graded free resolution $\mathbb{F}(I)$, which is unique up to isomorphism and has finite length at most $n$. Explicitly, $\mathbb{F}(I)$ can be written as
\begin{small}
\[
0 \rightarrow \bigoplus_{j \in \mathbb{Z}} R(-j)^{\beta_{\ell,j}} \xrightarrow{d_\ell} \cdots 
\rightarrow \bigoplus_{j \in \mathbb{Z}} R(-j)^{\beta_{i,j}} \xrightarrow{d_i} \cdots 
\rightarrow \bigoplus_{j \in \mathbb{Z}} R(-j)^{\beta_{0,j}} \xrightarrow{d_0} I \rightarrow 0,
\]
\end{small}
where $\ell \leq n$. The integers $\beta_{i,j}$ are the \emph{graded Betti numbers} of $I$.  
The \emph{Castelnuovo--Mumford regularity} (or simply \emph{regularity}) of $I$ is defined by
\[
\mathrm{reg}(I) = \max\{\, j \mid \beta_{i,\, i+j} \neq 0 \text{ for some } i \,\}.
\]
Moreover, one has $\mathrm{reg}(I) = \mathrm{reg}(R/I) + 1$. The quotient $R/I$ inherits a natural grading as a $K$-algebra, $R/I = \bigoplus_{k \in \mathbb{N}} (R/I)_k.$ The associated formal power series $\mathrm{HP}_{R/I}(t) = \sum_{k \in \mathbb{N}} \dim_K (R/I)_k\, t^k$ is called the \emph{Hilbert--Poincar\'e series} of $R/I$.  
By the Hilbert--Serre Theorem, there exists a unique polynomial $h(t) \in \mathbb{Z}[t]$, called the \emph{$h$-polynomial} of $R/I$, such that $h(1) \neq 0$ and
\[
\mathrm{HP}_{R/I}(t) = \frac{h(t)}{(1 - t)^d},
\]
where $d$ is the Krull dimension of $R/I$. Furthermore, if $R/I$ is Cohen--Macaulay, then $\mathrm{reg}(R/I) = \deg h(t)$.

We now introduce a variant of the standard rook polynomial, which we simply refer to as the \textit{rook polynomial} of $\cQ$.

Two rooks $R_1$ and $R_2$ are said to be in \textit{attacking position}, or are \textit{attacking rooks} in $\cQ$, if there exist two cells $A_1$ and $A_2$ of $\cQ$ in horizontal or vertical alignment containing $R_1$ and $R_2$, respectively, and such that the segment $[A_1,A_2]$ is contained in $\cQ$. Conversely, two rooks are in \textit{non-attacking position}, or are \textit{non-attacking rooks} in $\cQ$, if they are not in attacking position. For instance, in Figures~\ref{Figure: exa attacking rooks}(A) and (B), both configurations illustrate two non-attacking rooks. The function  \texttt{isNonAttackingRooks(A, B, Q)} checks whether two given rooks, encoded by the lower–left corners of the cells in which they are placed, are in a non-attacking position.

A \textit{$j$-rook configuration} in $\cQ$ is a set of $j$ rooks placed in non-attacking positions within $\cQ$, where $j \geq 0$; by convention, the $0$-rook configuration is $\emptyset$. The \textit{rook number} $r(\cQ)$ is defined as the maximum number of rooks that can be placed in $\cQ$ in non-attacking positions. We denote by $\mathcal{R}(\cQ,k)$ the set of all $k$-rook configurations in $\cQ$, and set $r_k = \vert \mathcal{R}(\cQ,k) \vert$, for all $k \in \{0, \dots, r(\cQ)\}$ (with the convention $r_0 = 1$). The \textit{rook polynomial} of $\cQ$ is then the polynomial in $\mathbb{Z}_{\geq 0}[t]$ defined by
\[
r_{\cQ}(t) = \sum_{k=0}^{r(\cQ)} r_k t^k.
\]
For instance, the polyomino in Figure~\ref{Figure: exa attacking rooks} has
$r_{\cQ}(t) = t^3 + 5 t^2 + 5 t + 1$ and $r(\cQ) = 3$. 
 The functions \texttt{allNonAttackingRookConfigurations Q}, \texttt{rookNumber Q}, \texttt{rookPolynomial Q} compute all non-attacking rook configurations in \texttt{Q}, the rook number and the rook polynomial of \texttt{Q}.

\begin{figure}[h]
\centering
\includegraphics[width=0.4\columnwidth]{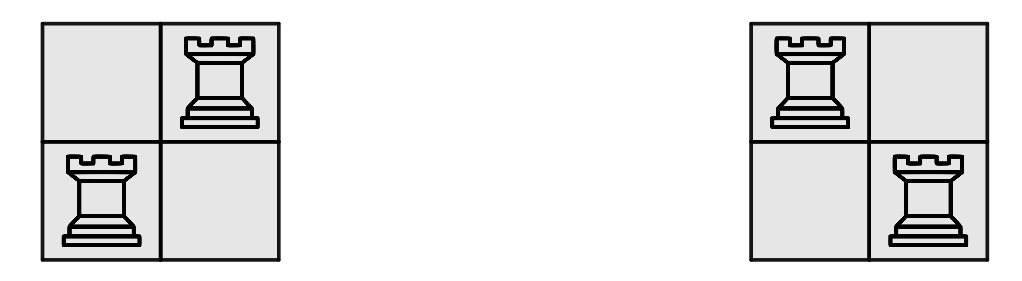}
\caption{Square tetromino and two equivalent $2$-rook configurations.}
\label{Fig: square}
\end{figure}

For some classes of collections of cells that do not contain the square tetromino (see Figure~\ref{Fig: square}), also known as \textit{thin}, the $h$-polynomial of the coordinate ring coincides with the rook polynomial of the associated collection of cells; see \cite{CJN, CNU2, DN2, N, RR}.

\subsection*{The switching rook polynomial.} In the non-thin case, the $h$-polynomial of $K[\cQ]$ differs from the rook polynomial of $\cQ$. Indeed, for a square tetromino $\cS$, we have $r_{\mathcal{S}}(t) = 1 + 4t + 2t^2$ and $h_{K[\mathcal{S}]}(t) = 1 + 4t + t^2$. In such cases, it coincides with a variant of the rook polynomial, called the \textit{switching rook polynomial}, which we now introduce. This connection has been established in several works, as \cite{EHQR, JN, KV, NQR1, NQR2, QRR}.

Let us note that $\bigcup_{j=0}^{r(\cQ)} \cR_j(\cQ)$ forms a simplicial complex, known as the \textit{chessboard complex} of $\cQ$. Two non-attacking rooks in $\cQ$ are said to be in \textit{switching position}, or are \textit{switching rooks}, if they occupy cells that are diagonally (or anti-diagonally) opposite within an \textit{inner interval} $I$ of $\cQ$, denoted $\cQ_I$. In this case, we say that the rooks are in a \textit{diagonal} (or \textit{anti-diagonal}) position. 

Fix $j \in \{0, \dots, r(\cQ)\}$. Let $F \in \cR_j(\cQ)$. Consider two switching rooks, $R_1$ and $R_2$, within $F$, positioned diagonally (or anti-diagonally) in $\cQ_I$ for some inner interval $I$. Let $R_1'$ and $R_2'$ be the rooks in the anti-diagonal (or diagonal, respectively) cells of $\cQ_I$. Then the set $(F \backslash \{R_1, R_2\}) \cup \{R_1', R_2'\}$ also belongs to $\cR_j(\cQ)$. This operation of replacing $R_1$ and $R_2$ with $R_1'$ and $R_2'$ is called a \textit{switch of $R_1$ and $R_2$}. 

\begin{figure}
\hfill \begin{minipage}[t]{0.22\linewidth}
\includegraphics[width=0.69\linewidth]{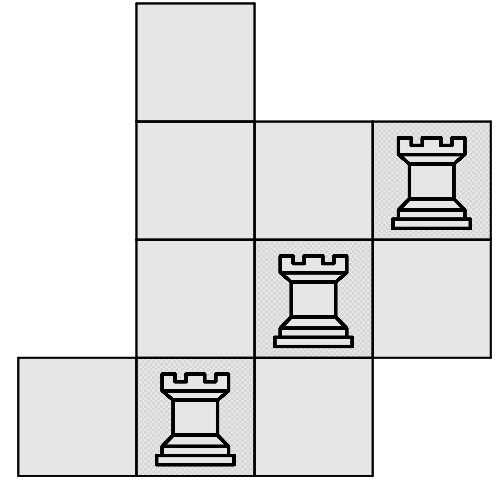}
\end{minipage}
\begin{minipage}[t]{0.22\linewidth}
\includegraphics[width=0.69\linewidth]{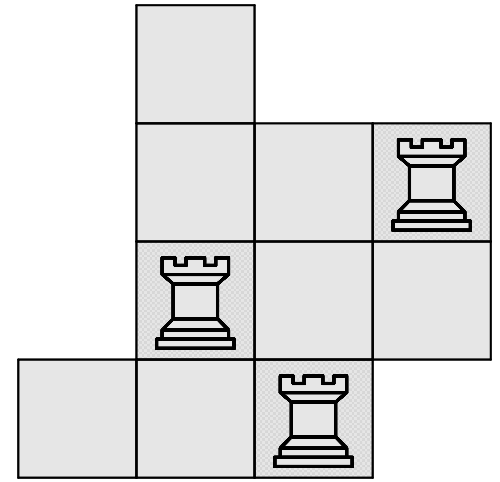}
\end{minipage}
\begin{minipage}[t]{0.22\linewidth}
\includegraphics[width=0.69\linewidth]{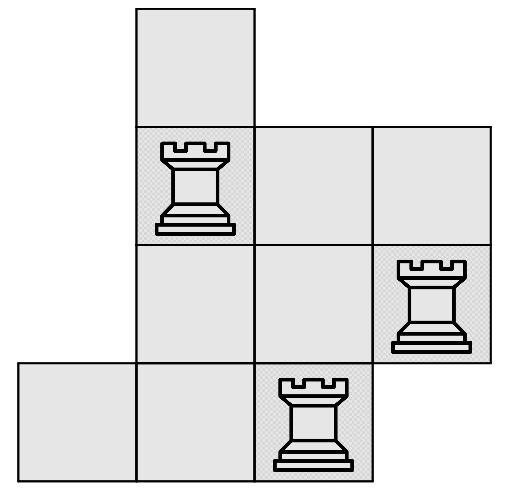}
\end{minipage}
\begin{minipage}[t]{0.22\linewidth}
\includegraphics[width=0.69\linewidth]{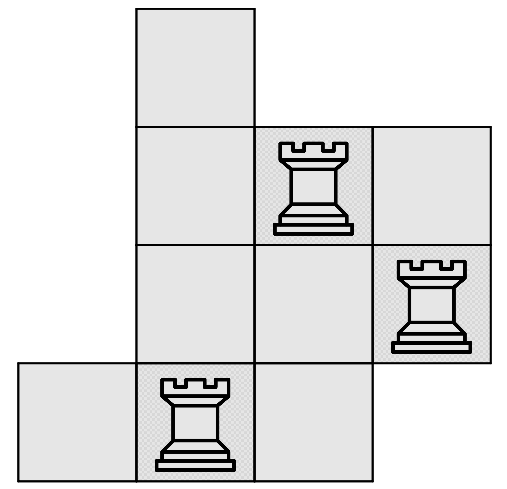}
\end{minipage}
\caption{Four arrangements of $3$ non-attacking rooks, equivalent under $\sim$.}
\label{Figure: equivalent}
\end{figure}

This defines an equivalence relation $\sim$ on $\cR_j(\cQ)$: we write $F_1 \sim F_2$ if $F_2$ can be obtained from $F_1$ through a sequence of switches. In this case, we say that $F_1$ and $F_2$ are \textit{equivalent with respect to $\sim$} (or \textit{same up to switches}). In Figure~\ref{Figure: equivalent}, four $3$-rook configurations equivalent under $\sim$ are shown.

Let $\tilde{\cR}_j(\cQ) = \cR_j(\cQ)/\sim$ be the set of equivalence classes. We define $\tilde{r}_j(\cQ) = \vert \tilde{\cR}_j(\cQ) \vert$ for $j \in \{0, \dots, r(\cQ)\}$, with the convention $\tilde{r}_0(\cQ) = 1$. The \textit{switching rook polynomial} of $\cQ$ is then defined as the polynomial in $\mathbb{Z}_{\geq 0}[t]$:
\[
\tilde{r}_{\cQ}(t) = \sum_{j=0}^{r(\cQ)} \tilde{r}_j(\cQ) t^j.
\]

  The functions \texttt{equivalenceClassesSwitchingRook Q} and \texttt{switchingRookPolynomial Q} compute the equivalence classes of non-attacking rook configurations, under switching, of \texttt{Q} and the switching rook polynomial of \texttt{Q}.
  



\subsection*{An example.} We conclude the section with the following example:

\begin{footnotesize}
\begin{lstlisting}
i1 : Q = cellCollection {{1,1}, {1,2}, {2,1}, {3,1}, {3,2}}
i2 : standardRookNumber Q
o2 = 2
i3 : standardRookPolynomial Q
       2
o3 = 4t  + 5t + 1
i4 : rookPolynomial Q
      3     2
o4 = t  + 5t  + 5t + 1
i5 : rookNumber Q
o5 = 3
i6 : Q = cellCollection {{1,1}, {1,2}, {2,1}, {2,2}}
i7 : switchingRookPolynomial Q
      2
o7 = t  + 4t + 1

i8 : rookPolynomial Q
      2
o8 = 2t  + 4t + 1
\end{lstlisting}
\end{footnotesize}

\end{document}